\documentclass[graybox]{svmult}

\usepackage{type1cm}        
\usepackage{makeidx}         
\usepackage{graphicx}        
\usepackage{multicol}        
\usepackage[bottom]{footmisc}

\usepackage{newtxtext}       %
\usepackage{newtxmath}       

\makeindex             

\begin{document}
\allowdisplaybreaks[4]
\title*{Finite Element Methods for  Elliptic Distributed Optimal Control Problems with Pointwise
 State Constraints}
\titlerunning{Finite Element Methods for  Elliptic Distributed Optimal Control Problems}
\author{Susanne C. Brenner}
\institute{Susanne C. Brenner \at Department of Mathematics and Center for Computation \& Technology,
Louisiana State University, Baton Rouge, LA 70803, USA. \email{brenner@math.lsu.edu}
}
%
%
\maketitle
\abstract*{Finite element methods for a model  elliptic distributed
 optimal control problem with pointwise state constraints are considered
 from the perspective of fourth order boundary value problems.}

\abstract{Finite element methods for a model  elliptic distributed
 optimal control problem with pointwise state constraints are considered
 from the perspective of fourth order boundary value problems.}
\section{Model Problem}
\label{sec:Model Problem}
 Let $\Omega$ be a convex bounded polygonal/polyhedral domain in
 $\mathbb{R}^2/\mathbb{R}^3$,
  $y_d\in L_2(\Omega)$, $\beta$ be a positive constant,
  $\psi\in H^3(\Omega)\cap W^{2,\infty}(\Omega)$ and $\psi>0$ on $\partial\Omega$.
    The model problem  \cite{Casas:1986:Control} is to find
\begin{equation}\label{brenner:OCP}
 (\bar y,\bar u)=\mathop{\rm argmin}_{(y,u)\in \mathbb{K}}
 \frac12\Big[\|y-y_d\|_{L_2(\Omega)}^2
   +\beta\|u\|_{L_2(\Omega)}^2\Big],
\end{equation}
 where $(y,u)\in H^1_0(\Omega)\times L_2(\Omega)$ belongs to $\mathbb{K}$ if and only if
\begin{alignat}{3}
 \int_\Omega\nabla y\cdot\nabla z\,dx&=
 \int_\Omega uz\,dx&\qquad&\forall\,z\in H^1_0(\Omega),
 \label{brenner:PDEConstraint}\\
 y&\leq\psi &\qquad&\text{a.e. on $\Omega$}.
 \label{brenner:StateConstraint}
\end{alignat}
 Throughout this paper we will follow the standard notation for operators,
 function spaces and norms
 that can be found for example in \cite{Ciarlet:1978:FEM,BScott:2008:FEM}.
\par
 In this model problem $y$ (resp., $u$) is the state (resp., control) variable,
 $y_d$ is the  desired state and $\beta$ is a regularization parameter.
  Similar linear-quadratic optimization problems also appear as subproblems
 when general PDE constrained optimization problems are solved by
  sequential quadratic programming (cf. \cite{HPUU:2009:Book,Troltzsch:2010:OC}).
\par
 In view of the convexity of $\Omega$, the constraint
 \eqref{brenner:PDEConstraint} implies $y\in H^2(\Omega)$
 (cf. \cite{Grisvard:1985:EPN,Dauge:1988:EBV,MR:2010:Polyhedral}).
 Therefore we can reformulate \eqref{brenner:OCP}--\eqref{brenner:StateConstraint}
 as follows:
\begin{equation}\label{brenner:OP}
 \text{Find}\quad\bar y=\mathop{\rm argmin}_{y\in K}\frac12
 \Big[\|y-y_d\|_{L_2(\Omega)}^2
   +\beta\|\Delta y\|_{L_2(\Omega)}^2\Big],
\end{equation}
 where
\begin{equation}\label{brenner:KDef}
 K=\{y\in H^2(\Omega)\cap H^1_0(\Omega): y\leq\psi \;\text{on $\Omega$}\}.
\end{equation}
\par
 Note that $K$ is nonempty because $\psi>0$ on $\partial\Omega$.
 It follows from the classical theory  of calculus of variations
 \cite{KS:1980:VarInequalities}
 that \eqref{brenner:OP}--\eqref{brenner:KDef} has a unique solution
 $\bar y\in K$
  characterized by the fourth order variational inequality
\begin{equation}\label{brenner:VI}
  a(\bar y,y-\bar y)\geq \int_\Omega y_d(y-\bar y)dx \qquad\forall\,y\in K,
\end{equation}
 where
\begin{equation}\label{brenenr:BilinearForm}
 a(y,z)=\beta\int_\Omega (\Delta y)(\Delta z)dx+\int_\Omega yz\,dx.
\end{equation}
 Furthermore, by the Riesz-Schwartz Theorem for nonnegative linear functionals
 \cite{Rudin:1966:RC,Schwartz:1966:Distributions},
 we can rewrite  \eqref{brenner:VI} as
\begin{equation}\label{brenner:Riesz}
  a(\bar y,z)=\int_\Omega y_dz\,dx+\int_\Omega z\,d\mu
  \qquad\forall\,z\in H^2(\Omega)\cap H^1_0(\Omega),
\end{equation}
 where
\begin{equation}\label{brenner:Multiplier}
 \text{$\mu$ is a nonpositive finite Borel measure}
\end{equation}
  that satisfies the complementarity condition
\begin{equation}\label{brenner:Complementarity}
 \int_\Omega (\bar y-\psi)d\mu=0.
\end{equation}
\par
 Note that \eqref{brenner:Complementarity} is equivalent to the statement that
\begin{equation}\label{brenner:Active}
  \text{$\mu$ is supported on $\mathcal{A}$},
\end{equation}
 where the active set $\mathcal{A}=\{x\in\Omega:\,\bar y(x)=\psi(x)\}$ satisfies
\begin{equation}\label{brenner:Compactness}
 \mathcal{A}\subset\subset \Omega
\end{equation}
 because $\psi>0$ on $\partial\Omega$ and $\bar y=0$ on $\partial\Omega$.
\par
 According to the elliptic regularity theory in
  \cite{Grisvard:1985:EPN,Dauge:1988:EBV,MR:2010:Polyhedral,Frehse:1971:VarInequality,Frehse:1973:VI},
  we have
\begin{equation}\label{brenner:StateRegularity}
 \bar y\in H^3_{loc}(\Omega)\cap W^{2,\infty}_{loc}(\Omega)\cap H^{2+\alpha}(\Omega),
\end{equation}
 where $\alpha\in (0,1]$ is determined by the geometry of $\Omega$.
 It then follows from  \eqref{brenner:Riesz},
  \eqref{brenner:Active}--\eqref{brenner:StateRegularity} and
  integration by parts that
\begin{equation}\label{brenner:MultiplierRegularity}
 \mu\in H^{-1}(\Omega).
\end{equation}
 Details for \eqref{brenner:StateRegularity} and \eqref{brenner:MultiplierRegularity}
 can be found in \cite{BGS:2018:POne}.
\begin{remark}\label{brenner:Alternative}
  Note that (cf. \cite{Ladyzhenskaya:1958:Elliptic,Grisvard:1985:EPN})
\begin{equation*}
  \int_\Omega (\Delta y)(\Delta z)dx=\int_\Omega D^2y:D^2 z\,dx\qquad
   \forall\,y,z\in H^2(\Omega)\cap H^1_0(\Omega),
\end{equation*}
 where $D^2 y:D^2 z$ denotes the Frobenius inner product between the Hessian matrices
 of $y$ and $z$.  Therefore we can rewrite the bilinear form $a(\cdot,\cdot)$
 in \eqref{brenenr:BilinearForm} as
\begin{equation}\label{brenner:AlternativeBilinearForm}
  a(y,z)= \beta\int_\Omega D^2y:D^2z\,dx+\int_\Omega yz\,dx.
\end{equation}
\end{remark}
\section{Finite Element Methods}
\label{brenner:FEMs}
 In the absence of the state constraint \eqref{brenner:StateConstraint}, we have
 $K=H^2(\Omega)\cap H^1_0(\Omega)$ and \eqref{brenner:VI} becomes the boundary value problem
\begin{equation}\label{brenner:BVP}
 a(\bar y,z)=\int_\Omega y_dz\,dx\qquad\forall\,z\in H^2(\Omega)\cap H^1_0(\Omega).
\end{equation}
 Since \eqref{brenner:BVP} is essentially a bending problem for simply supported plates,
 it can be solved by many finite element methods such as
 (i) conforming methods, (ii) classical nonconforming methods,
 (iii) discontinuous Galerkin methods, and (iv) mixed methods.
 For the sake of brevity, below we will consider these methods for
 $\Omega\subset \mathbb{R}^2$.
 But all the results can be extended to three dimensions.
\par
 Let $V_h$ be a finite element space associated with a triangulation $\mathcal{T}_h$
 of $\Omega$. The approximate solution
 $\bar y_h\in V_h$ is determined by
\begin{equation}\label{brenner:BVPFEM}
 a_h(\bar y_h,z)=\int_\Omega y_dz\,dx\qquad\forall\,z\in V_h,
\end{equation}
 where the choice of the bilinear form $a_h(\cdot,\cdot)$ depends
  on the type of finite element method being used.
\par\medskip\noindent
{\bf Conforming Methods}
\par\smallskip
 In this  case $V_h\subset H^2(\Omega)\cap H^1_0(\Omega)$ is a $C^1$ finite element space
 and we can take $a_h(\cdot,\cdot)$ to be $a(\cdot,\cdot)$.
 This class of methods includes
 the Bogner-Fox-Schmit element \cite{BFS:1965:Element},
 the Argyris elements \cite{AFS:1968:TUBA}, the macro elements
 \cite{CT:1965:Element,Ciarlet:1974:HCT,DDPS:1979:Macro}, and generalized finite elements
 \cite{MB:1996:PTU,BBO:2003:GFEM,ODJ:2011:Plate}.
\par\medskip\noindent{\bf Classical Nonconforming Methods}
\par\smallskip
  In this case $V_h\subset L_2(\Omega)$ consists of
  finite element functions that are weakly continuous up to
  first order derivatives across element boundaries,
  and the bilinear form $a_h(\cdot,\cdot)$ is given by
\begin{equation}\label{brenner:NCah}
 a_h(y,z)=\beta\sum_{T\in\mathcal{T}_h}\int_\Omega D^2y:D^2z\,dx
  +\int_\Omega yz\,dx.
\end{equation}
 Here we are using the piecewise version of \eqref{brenner:AlternativeBilinearForm}, which
 provides better local control of the nonconforming energy norm
 $\|\cdot\|_{a_h}=\sqrt{a_h(\cdot,\cdot)}$.
\par
 This class of methods includes the Adini element \cite{AC:1961:Plate},
 the Zienkiewicz element \cite{BCIZ:1965:TEB}, the Morley element \cite{Morley:1968:plate},
 the Fraeijs de Veubeke element \cite{FdV:1974:PatchTest},
 and the incomplete biquadratic element \cite{Shi:1986:IBPlate}.
%
\par\medskip\noindent{\bf Discontinuous Galerkin Methods}
\par\smallskip
  In this case $V_h$ consists of functions that are totally discontinuous or only discontinuous
  in the normal derivatives across element boundaries,
  and stabilization terms are included in the bilinear form $a_h(\cdot,\cdot)$.
  The simplest choice is a Lagrange finite element space $V_h\subset H^1_0(\Omega)$, resulting in
  the $C^0$ interior penalty methods \cite{EGHLMT:2002:DG3D,BSung:2005:DG4,Brenner:2012:C0IP},
  where the bilinear form $a_h(\cdot,\cdot)$ is given by
\begin{align}\label{brenner:C0IPah}
 a_h(y,z)& =\beta\bigg[\sum_{T\in\mathcal{T}_h}\int_T D^2y:D^2z\,dx
  +\sum_{e\in\mathcal{E}_h^i}
  \int_e\{\hspace{-2pt}\{\partial^2 y/\partial n^2\}\hspace{-2pt}\}
  [\hspace{-2pt}[\partial z/\partial n]\hspace{-2pt}]\,ds\notag\\
   &\hspace{20pt}+\sum_{e\in\mathcal{E}_h^i}
   \int_e{\{\hspace{-2pt}\{\partial^2 z/\partial n^2\}\hspace{-2pt}\}}
  [\hspace{-2pt}[\partial y/\partial n]\hspace{-2pt}]\,ds\\
  &\hspace{40pt}+\sigma\sum_{e\in\mathcal{E}_h^i} |e|^{-1}
  \int_e[\hspace{-2pt}[\partial y/\partial n]\hspace{-2pt}]
   [\hspace{-2pt}[\partial z/\partial n]\hspace{-2pt}]\,ds\bigg]
     +\int_\Omega yz\,dx.\notag
\end{align}
 Here $\mathcal{E}_h^i$ is the set of the interior edges of $\mathcal{T}_h$,
 $\{\hspace{-2pt}\{\partial^2 y/\partial n^2\}\hspace{-2pt}\}$
 (resp., $[\hspace{-2pt}[\partial y/\partial n]\hspace{-2pt}] $) is
 the average (resp., jump) of the
 second (resp., first) normal derivative of $y$ across the edge $e$,
 $|e|$ is the length of the edge $e$,
 and $\sigma$ is a (sufficiently large) penalty parameter.
\par
 Other discontinuous Galerkin methods for fourth order problems can be found in
 \cite{SM:2007:Biharmonic,HHH:2010:C0DG,HH:2016:C0IP}.
\par\medskip\noindent
{\bf Mixed Methods}
\par\smallskip
 In this case $V_h\subset H^1_0(\Omega)$ is a Lagrange finite element space.
  The approximate solution
 $\bar y_h$ is determined by
\begin{alignat}{3}
  \int_\Omega \bar y_hz\,dx+\beta\int_\Omega \nabla \bar u_h\cdot\nabla z\,dx&
  =\int_\Omega y_d z\,dx
  &\qquad&\forall\,z\in V_h,\label{brenner:Mixed1}\\
  \int_\Omega\nabla \bar y_h\cdot\nabla v\,dx-\int_\Omega \bar u_h v\,dx&=0
  &\qquad&\forall\,v\in V_h.\label{brenner:Mixed2}
\end{alignat}
 By eliminating $\bar u_h$ from \eqref{brenner:Mixed1}--\eqref{brenner:Mixed2},
 we can recast $\bar y_h$ as the solution of \eqref{brenner:BVPFEM} where
\begin{equation}\label{brenner:MixedAh}
 a_h(y,z)=\beta\int_\Omega (\Delta_h y)(\Delta_h z)\,dx+\int_\Omega yz\,dx,
\end{equation}
 and the discrete Laplace operator $\Delta_h:V_h\longrightarrow V_h$ is defined by
\begin{equation}\label{brenner:DiscreteLaplace}
   \int_\Omega (\Delta_h y)z\,dx=-\int_\Omega\nabla y\cdot\nabla z\,dx
   \qquad\forall\,y,z\in V_h.
\end{equation}
\par\medskip\noindent
{\bf Finite Element Methods for the Optimal Control Problem}
\par\smallskip
 With the finite element methods for \eqref{brenner:BVP} in hand,
 we can now simply discretize the variational inequality \eqref{brenner:VI} as
 follows: Find $\bar y_h\in V_h$
 such that
\begin{equation}\label{brenner:DVI}
 a_h(\bar y_h,y-\bar y_h)\geq \int_\Omega y_d(y-\bar y_h)dx\qquad\forall\,y\in K_h,
\end{equation}
 where
\begin{equation}\label{brenner:Kh}
 K_h=\{y\in V_h:\,I_h y\leq I_h \psi\;\text{on $\Omega$}\},
\end{equation}
 and $I_h$ is the nodal interpolation operator for the conforming $P_1$ finite element space
  associated with $\mathcal{T}_h$.
 In other words, the constraint \eqref{brenner:StateConstraint}
 is only imposed at the vertices of $\mathcal{T}_h$.
\begin{remark}\label{brenner:Implementations}
 Conforming, nonconforming, $C^0$ interior penalty and mixed methods for \eqref{brenner:VI}
 were investigated in
  \cite{BGS:2018:POne,Meyer:2008:OptimalControl,LGY:2009:Control,GY:2011:State,
  BSZ:2013:OptimalControl,CMV:2014:State,BDS:2014:PUMOC,BOPPSS:2016:OC3D}.
\end{remark}
\section{Convergence Analysis}\label{brenner:Convergence}
 For simplicity, we will only provide details for the case of conforming finite element methods
 and briefly describe the extensions to other methods at the end of the section.
\par
 For conforming finite element methods,
 we have $a_h(\cdot,\cdot)=a(\cdot,\cdot)$ and the energy norm $\|\cdot\|_a=\sqrt{a(\cdot,\cdot)}$
 satisfies  , by a Poincar\'e-Friedrichs inequality \cite{Necas:2012:Direct},
\begin{equation}\label{brenner:NormEquivalency}
 \|v\|_a\approx \|v\|_{H^2(\Omega)} \qquad \forall\,v\in H^2(\Omega).
\end{equation}
 Our goal is to show that
\begin{equation}\label{brenner:ErrorEstimate}
  \|\bar y-\bar y_h\|_a\leq Ch^\alpha,
\end{equation}
 where $\alpha$ is the index of elliptic regularity that appears in
 \eqref{brenner:StateRegularity}.
\par
  We assume (cf. \cite{GS:2002:Interpolation}) that
 there exists an operator $\Pi_h:H^2(\Omega)\cap H^1_0(\Omega)\longrightarrow V_h$ such that
\begin{equation}\label{brenner:PihNodal}
 \Pi_h\zeta=\zeta \quad\text{at the vertices of $\mathcal{T}_h$}
\end{equation}
 and
\begin{equation}\label{brenner:PihEstimate}
   \|\zeta-\Pi_h\zeta\|_{L_2(\Omega)}+h|\zeta-\Pi_h\zeta|_{H^1(\Omega)}
   +h^2|\zeta-\Pi_h\zeta|_{H^2(\Omega)}\leq Ch^{2+\alpha}|\zeta|_{H^{2+\alpha}(\Omega)}
\end{equation}
 for all $\zeta\in H^{2+\alpha}(\Omega)\cap H^1_0(\Omega)$,
 where $h=\max_{T\in\mathcal{T}_h}\text{diam}\,T$ is the mesh size of the triangulation
 $\mathcal{T}_h$.  Here and below we use $C$ to denote a generic positive constant independent of
 $h$.
\par
 In particular \eqref{brenner:KDef}, \eqref{brenner:Kh} and \eqref{brenner:PihNodal} imply
\begin{equation}\label{eq:Mapping}
  \text{$\Pi_h$ maps $K$ into $K_h$.}
\end{equation}
\goodbreak
 Therefore $K_h$ is nonempty and the discrete problem defined by
 \eqref{brenner:DVI}--\eqref{brenner:Kh} has a unique solution.
\par
 We will also use the following standard properties of the interpolation operator
 $I_h$ (cf. \cite{Ciarlet:1978:FEM,BScott:2008:FEM}):
\begin{alignat}{3}
  \|\zeta- I_h\zeta\|_{L_\infty(T)}&\leq Ch_T^2|\zeta|_{W^{2,\infty}(T)}&\qquad&
  \forall\,\zeta\in W^{2,\infty}(T),\, T\in \mathcal{T}_h, \label{brenner:Infinity}\\
  |\zeta-I_h\zeta|_{H^1(T)}&\leq Ch_T|\zeta|_{H^2(T)}&\qquad&\forall\,\zeta\in H^2(T),\,
   T\in \mathcal{T}_h,
  \label{brenner:HOne}
\end{alignat}
 where $h_T$ is the diameter of $T$.
\par
 We begin with the estimate
\begin{align}\label{brenner:Preliminary}
 \|\bar y-\bar y_h\|_a^2&=a(\bar y-\bar y_h,\bar y-\bar y_h)\notag \\
   &=a(\bar y-\bar y_h,\bar y-\Pi_h\bar y)+a(\bar y,\Pi_h\bar y-\bar y_h)-
      a(\bar y_h,\Pi_h\bar y-\bar y_h)\\
         &\leq C_1\|\bar y-\bar y_h\|_ah^\alpha+\Big[a(\bar y,\Pi_h\bar y-\bar y_h)
                 -\int_\Omega y_d(\Pi_h\bar y-\bar y_h)dx\Big]\notag
\end{align}
 that follows from \eqref{brenner:StateRegularity}, \eqref{brenner:DVI},
 \eqref{brenner:NormEquivalency}, \eqref{brenner:PihEstimate}, \eqref{eq:Mapping} and
 the Cauchy-Schwarz inequality.
\begin{remark}\label{brenner:Analogous}
 Note that an estimate analogous to \eqref{brenner:Preliminary} also
 appears in the error analysis
 for the boundary value problem \eqref{brenner:BVP}.  Indeed the second term
 on the right-hand side
 of \eqref{brenner:Preliminary} vanishes in the case of  \eqref{brenner:BVP}
 and we would have
 arrived at the desired estimate $\|\bar y-\bar y_h\|_a\leq Ch^\alpha$.
\end{remark}
\par
 The idea now is to show that
\begin{equation}\label{brenner:KeyEstimate}
 a(\bar y,\Pi_h\bar y-\bar y_h)
      -\int_\Omega y_d(\Pi_h\bar y-\bar y_h)dx
       \leq C_2\big[h^{2\alpha}+h^\alpha\|\bar y-\bar y_h\|_a\big],
\end{equation}
 which together with \eqref{brenner:Preliminary} implies
\begin{equation}\label{brenner:Young}
 \|\bar y-\bar y_h\|_a^2\leq C_3h^\alpha\|\bar y-\bar y_h\|_a+C_2 h^{2\alpha}.
\end{equation}
 The estimate \eqref{brenner:ErrorEstimate} then follows from
 \eqref{brenner:Young} and the inequality
\begin{equation*}
  ab\leq \frac{\epsilon}{2}a^2+\frac{1}{2\epsilon}b^2
\end{equation*}
 that holds for any positive $\epsilon$.
\par
 Let us turn to the derivation of \eqref{brenner:KeyEstimate}.
 Since $K_h\subset  V_h\subset H^2(\Omega)\cap H^1_0(\Omega)$, we have,
 according to \eqref{brenner:Riesz},
\begin{align}\label{brenner:K1}
 & a(\bar y,\Pi_h\bar y-\bar y_h)
                          -\int_\Omega y_d(\Pi_h\bar y-\bar y_h)dx
         =\int_\Omega (\Pi_h\bar y-\bar y_h)d\mu \notag \\
    &\hspace{70pt}=\int_\Omega (\Pi_h\bar y-\bar y)d\mu+\int_\Omega (\bar y-\psi)d\mu
           +\int_\Omega (\psi-I_h\psi)d\mu\\
          &\hspace{100pt}  +\int_\Omega (I_h\psi-I_h\bar y_h)d\mu  +
           \int_\Omega (I_h\bar y_h-\bar y_h)d\mu,   \notag
\end{align}
\goodbreak
 and, in view of \eqref{brenner:Multiplier},
 \eqref{brenner:Complementarity} and \eqref{brenner:Kh},
\begin{equation}\label{brenner:K2}
  \int_\Omega (\bar y-\psi)d\mu=0 \quad \text{and} \quad
   \int_\Omega (I_h\psi-I_h\bar y_h)d\mu\leq 0.
\end{equation}
\par
 We can estimate the other three integrals on the right-hand side of
  \eqref{brenner:K1} as follows:
\begin{equation}\label{brenner:K3}
  \int_\Omega (\Pi_h\bar y-\bar y)d\mu\leq \|\mu\|_{H^{-1}(\Omega)}
    \|\Pi_h\bar y-\bar y\|_{H^1(\Omega)}\leq Ch^{1+\alpha}
\end{equation}
 by \eqref{brenner:StateRegularity}, \eqref{brenner:MultiplierRegularity} and
 \eqref{brenner:PihEstimate};
\begin{equation}\label{brenner:K4}
 \int_\Omega (\psi-I_h\psi)d\mu\leq |\mu(\Omega)|\|\psi-I_h\psi\|_{L_\infty(\Omega)}\\
   \leq Ch^2
\end{equation}
 by \eqref{brenner:Multiplier} and \eqref{brenner:Infinity};
\begin{align}\label{brenner:K5}
   &\int_\Omega (I_h\bar y_h-\bar y_h)d\mu=
      \int_\Omega \big[I_h(\bar y_h-\bar y)-(\bar y_h-\bar y)\big]d\mu+
         \int_\Omega (I_h\bar y-\bar y)d\mu\notag\\
      &\hspace{20pt}\leq \|\mu\|_{H^{-1}(\Omega)}
       |I_h(\bar y_h-\bar y)-(\bar y_h-\bar y)|_{H^1(\Omega)}
            +|\mu(\Omega)| \|I_h\bar y-\bar y\|_{L_\infty(\mathcal{A})}\\
            &\hspace{20pt}\leq C\big[h|\bar y_h-\bar y|_{H^2(\Omega)}+ h^2\big]\notag\\
           &\hspace{20pt} \leq C\big(h\|\bar y-\bar y_h\|_a+h^2\big]  \notag
\end{align}
 by \eqref{brenner:Active}--\eqref{brenner:StateRegularity},
 \eqref{brenner:NormEquivalency}, \eqref{brenner:Infinity} and \eqref{brenner:HOne}.
\par
 The estimate \eqref{brenner:KeyEstimate} follows from \eqref{brenner:K1}--\eqref{brenner:K5}
 and the fact that $\alpha\leq 1$.
\par\medskip
 The estimate \eqref{brenner:ErrorEstimate} can be extended to the other
 finite element methods in Section~\ref{brenner:FEMs} provided $\|\cdot\|_a$ is replaced by $\|\cdot\|_{a_h}=\sqrt{a_h(\cdot,\cdot)}$.
\par
  For classical nonconforming finite
 element methods and discontinuous Galerkin methods, the key ingredient
  for the convergence analysis, in addition to an
  operator $\Pi_H:H^2(\Omega)\cap H^1_0(\Omega)\longrightarrow V_h$ that satisfies
  \eqref{brenner:PihNodal} and \eqref{brenner:PihEstimate},   is the existence of
 an {\em enriching} operator $E_h:\longrightarrow H^2(\Omega)\cap H^1_0(\Omega)$
 with the following properties:
\begin{align}
  &(E_hv)(p)=v(p) \quad\text{for all vertices $p$ of $\mathcal{T}_h$},\label{brenner:Eh1}\\
   &\|v-E_hv\|_{L_2(\Omega)}+h\Big(\sum_{T\in\mathcal{T}_h}|v-E_hv|_{H^1(T)}^2\Big)^\frac12
   +h^2|E_hv|_{H^2(\Omega)}\notag\\
      &\hspace{140pt}\leq Ch^2 \|v\|_h\quad\forall\,v\in V_h,\label{brenner:Eh2}\\
  &\|\zeta-E_h\Pi_h\zeta\|_{H^1(\Omega)}\leq Ch^{1+\alpha}\|\zeta\|_{H^{2+\alpha}(\Omega)}
  \qquad\forall\,\zeta\in H^{2+\alpha}(\Omega)\cap H^1_0(\Omega),\label{brenner:Eh3}\\
  &|a_h(\Pi_h\zeta,v)-a(\zeta,E_hv)|\leq Ch^\alpha\|\zeta\|_{H^{2+\alpha}(\Omega)}\|v\|_h
  \label{brenner:Eh4}
\end{align}
 for all $\zeta\in H^{2+\alpha}(\Omega)\cap H^1_0(\Omega)$ and $v\in V_h$.
\par
 Property \eqref{brenner:Eh1} is related to the fact that the discrete constraints are
 imposed at the vertices of $\mathcal{T}_h$; property \eqref{brenner:Eh2} indicates that
 in some sense $\|v-E_hv\|_h$  measures the distance between $V_h$ and $H^2(\Omega)\cap H^1_0(\Omega)$;
 property \eqref{brenner:Eh3} means that $E_h\Pi_h$ behaves like a quasi-local interpolation operator;
 property \eqref{brenner:Eh4} states that $E_h$ is essentially
 the adjoint of $\Pi_h$ with respect to
 the continuous and discrete bilinear forms.  The idea is to use \eqref{brenner:Eh2} and
 \eqref{brenner:Eh4} to reduce the error estimate to the continuous level, and then the error
 analysis can proceed as in the case of conforming finite element method by using
 \eqref{brenner:Eh1} and \eqref{brenner:Eh3}.  Details can be found in
 \cite{BSung:2017:State}.
\begin{remark}\label{brenner:Enriching}
  The operator $E_h$ maps $V_h$ to a conforming finite element space
  and its construction is based on averaging.  The history of using such enriching
  operators to handle nonconforming finite element methods is discussed in
  \cite{Brenner:2015:CR}.
\end{remark}
\par
 In the case of the mixed method where $V_h\subset H^1_0(\Omega)$ is a Lagrange finite
 element space, the operator
  $E_h:V_h\longrightarrow H^2(\Omega)\cap H^1_0(\Omega)$
 is defined by
\begin{equation}\label{brenner:MixedEh}
 \int_\Omega \nabla E_hv\cdot\nabla w\,dx=\int_\Omega \nabla v\cdot\nabla w\,dx
 \qquad\forall v\in V_h,\,w\in H^1_0(\Omega).
\end{equation}
 The properties \eqref{brenner:Eh2}--\eqref{brenner:Eh4} remain valid provided
 $\Pi_h$ is replaced by the Ritz projection operator $R_h:H^1_0(\Omega)\longrightarrow V_h$
 defined by
\begin{equation}\label{brenner:Ritz}
  \int_\Omega\nabla R_h\zeta\cdot\nabla v\,dx=\int_\Omega\nabla\zeta\cdot\nabla v\,dx
  \qquad\forall\,v\in V_h.
\end{equation}
 In fact \eqref{brenner:MixedEh} and \eqref{brenner:Ritz} imply $\zeta-E_hR_h\zeta=0$  and
 property \eqref{brenner:Eh3} becomes trivial.  However the properties \eqref{brenner:PihNodal}
 and \eqref{brenner:Eh1} no longer hold, which necessitates the use of the more sophisticated
 interior error estimates (cf. \cite{Wahlbin:2000:Handbook})
 in the convergence analysis.
 Details can be found in \cite{BGS:2018:POne}.
\begin{remark}\label{brenner:Other}
  Since the elliptic regularity
  index $\alpha$ in \eqref{brenner:StateRegularity} is determined by the
   singularity of the Laplace equation
  near the boundary of $\Omega$, various finite element techniques
  \cite{FGW:1973:SFM,BKP:1979:DIE} can be employed
  to improve the estimate \eqref{brenner:ErrorEstimate} to
\begin{equation}\label{brenner:Improvement}
 \|\bar y-\bar y_h\|_{a_h}\leq Ch.
\end{equation}
 One can also compute an approximation $\bar u_h$ for the optimal control $\bar u$ from the
 approximate optimal state $\bar y_h$  through post-processing processes
 \cite{BSZ:2015:PP}.
\end{remark}
\begin{remark}\label{brenner:PDAS}
 The discrete problems generated by the finite element methods in
  Section~\ref{brenner:FEMs}, which only involve simple box constraints, can be
  solved efficiently by a primal-dual active set algorithm
  \cite{BK:2002:PDAS,HIK:2003:PDAS,IK:2008:Lagrange}.
\end{remark}
\section{Concluding Remarks}
\label{brenner:Conclusion}
 In this paper finite element methods for elliptic  distributed optimal control
 problems with pointwise state constraints are treated from the perspective of
 finite element methods for the boundary value problem of simply supported plates.
\par
 The discussion in Section~\ref{brenner:FEMs} shows that one can solve
 elliptic distributed optimal control problems with pointwise state constraints by a
 straight-forward adaptation of many finite element methods for simply supported
 plates.   The convergence analysis in
 Section~\ref{brenner:Convergence} demonstrates that the gap between the finite
 element analysis for boundary value problems and the finite element analysis for
 elliptic optimal control problems is in fact quite narrow.  Thus the vast arsenal of
 finite element techniques developed
  for elliptic boundary value problems over several decades can be
 applied to elliptic optimal control problems with only minor modifications.
\par
 Note that in the traditional approach to elliptic optimal control problems,
 the optimal control $\bar u$ is treated as the primary unknown and the resulting
 finite element methods in \cite{Meyer:2008:OptimalControl,CMV:2014:State}
 are equivalent to the method defined by
 \eqref{brenner:DVI}, where the bilinear form is given by
 \eqref{brenner:MixedAh}.
 Therefore the approach
 based on the reformulation \eqref{brenner:OP}--\eqref{brenner:KDef}
 expands the scope of finite element methods for elliptic optimal control
 problems from a special class of methods (i.e., mixed methods) to
 all classes of methods.
 In addition to the finite element mentioned in Section~\ref{brenner:FEMs}, one can also consider
  recently developed finite element methods for fourth order problems on polytopal meshes
  \cite{BM:2013:VEM,MWY:2014:WG4,WW:2014:WG4,CM:2016:VEM,AMV:2018:VEM4,
  ZZCM:2018:MorleyVEM,BDGK:2018:HHOPlate,BDR:2019:3DC1}.
\par
 The new approach has been extended to problems with the Neumann boundary condition
 \cite{BSZ:2019:Neumann,BOS:2019:Neumann}
 and to problems
 with pointwise constraints on both
 control and state \cite{BGPS:2018:Morley}.  It has also been extended
 to problems on nonconvex domains
 \cite{BGS:2018:POne,BGS:2018:Nonconvex,BOS:2019:Neumann}.
\par
  Below are some open problems related to the finite element methods presented in
  Section~\ref{brenner:FEMs}.
\begin{enumerate}
  \item It follows from the error estimates \eqref{brenner:ErrorEstimate} and
  \eqref{brenner:Improvement} that
\begin{equation}\label{brenner:LowerOrderNorm}
  \|\bar y-\bar y_h\|_{H^1(\Omega)}+\|\bar y-\bar y_h\|_{L_\infty(\Omega)}
   \leq Ch^\gamma,
\end{equation}
 where $\gamma=\alpha$ (without special treatment) or $1$ (with special treatments).
 For conforming or mixed finite element methods, the
 estimate \eqref{brenner:LowerOrderNorm} is a direct consequence of the
 fact that the energy norm is equivalent to the $H^2(\Omega)$ norm and that we have
 the Sobolev inequality
\begin{equation*}
  \|\zeta\|_{L_\infty(\Omega)}\leq C\|\zeta\|_{H^2(\Omega)}.
\end{equation*}
 For classical nonconforming and discontinuous Galerkin methods, the
 estimate \eqref{brenner:LowerOrderNorm} follows from the Poincar\'e-Friedrichs inequality
 and Sobolev inequality for piecewise $H^2$ functions in
 \cite{BWZ:2004:PF4,BNRS:2017:VonKarman}.
\par\hspace{12pt}
  Comparing to $\|\cdot\|_{H^2(\Omega)}$, the norms $\|\cdot\|_{H^1(\Omega)}$ and
  $\|\cdot\|_{L_\infty(\Omega)}$ are lower order norms and, based on
  experience with finite element methods for the boundary value problem \eqref{brenner:BVP},
  the convergence in $\|\cdot\|_{H^1(\Omega)}$ and
  $\|\cdot\|_{L_\infty(\Omega)}$ should be of higher order,
  and this is observed in numerical experiments.  But the theoretical justifications
  for the observed higher order convergence is missing.  In the case of the boundary value problem
  \eqref{brenner:BVP}, one can show higher order convergence for lower order norms through
  a duality argument.  However duality arguments do not work for variational inequalities
  even in one dimension \cite{CC:2018:Obstacle}.
  New ideas are needed.\\
\item An interesting phenomenon concerning fourth order variational inequalities is that
   {\em a posteriori} error estimators
 originally designed for fourth order boundary value problems can be
 directly applied to fourth order
 variational inequalities \cite{BSZ:2019:Neumann,BGSZ:2017:AdaptiveVI4}.
  This is
 different from the second order case where {\em a posteriori} error estimators for boundary
 value problems are not directly applicable to variational inequalities. This difference is
 essentially due to the fact that Dirac point measures belong to $H^{-2}(\Omega)$ but not
 $H^{-1}(\Omega)$.
\par\hspace{12pt}
 Optimal convergence of these adaptive finite element methods have been observed in numerical
 experiments.  However the proofs of convergence and optimality are missing.\\
\item Fast solvers  for fourth order variational
 inequalities is an almost completely open area.  Some recent work on additive Schwarz
 preconditioners for the subsystems that appear in the primal-dual active set algorithm
 can be found in \cite{BDS:2018:ASP,BDS:2018:DD25}.  Much remains to be done.
\end{enumerate}

\begin{acknowledgement}
This paper is based on research supported by the National Science Foundation under Grant Nos.
DMS-13-19172, DMS-16-20273 and DMS-19-13035.
\end{acknowledgement}

\begin{thebibliography}{99.}
%
\bibitem{Casas:1986:Control} Casas, E.: Control of an elliptic problem with pointwise state
              constraints.
   SIAM J. Control Optim. \textbf{24}, 1309--1318 (1986)
%
\bibitem{Ciarlet:1978:FEM}
  Ciarlet, P.G.:
  The Finite Element Method for Elliptic Problems.
  North-Holland,  Amsterdam (1978)
%
\bibitem{BScott:2008:FEM}
  Brenner, S.C.,  Scott, L.R.:
  The Mathematical Theory of Finite Element Methods (Third Edition).
  Springer-Verlag, New-York (2008)
%
\bibitem{HPUU:2009:Book}
    Hinze, M. and Pinnau, R. and Ulbrich, M. and Ulbrich, S.:
    Optimization with PDE Constraints.
   Springer, New York (2009)
%
\bibitem{Troltzsch:2010:OC}
    Tr{\"o}ltzsch, F.:
     Optimal Control of Partial Differential Equations.
  American Mathematical Society, Providence (2010)
%
\bibitem{Grisvard:1985:EPN}
 Grisvard, P.:
  Elliptic Problems in Non Smooth Domains.
 Pitman, Boston (1985)
%
\bibitem{Dauge:1988:EBV}
  Dauge, M.: Elliptic Boundary Value Problems on Corner
                   Domains.   Springer-Verlag, Berlin-Heidelberg (1988)
%
\bibitem{MR:2010:Polyhedral}
    Maz'ya, V. , Rossmann, J.:
     Elliptic Equations in Polyhedral Domains.
  American Mathematical Society,  Providence (2010)
%
\bibitem{KS:1980:VarInequalities}
    Kinderlehrer, D., Stampacchia, G.:
    An Introduction to Variational Inequalities and Their
              Applications.
  Society for Industrial and Applied Mathematics,
   Philadelphia (2000)
%
\bibitem{Rudin:1966:RC}
    Rudin, W.: Real and Complex Analysis.
   McGraw-Hill, New York (1966)
%
\bibitem{Schwartz:1966:Distributions}
    Schwartz, L.:
     Th\'eorie des Distributions. Hermann, Paris (1966)
%
\bibitem{Frehse:1971:VarInequality}
    Frehse, J.: Zum {D}ifferenzierbarkeitsproblem bei
              {V}ariationsungleichungen h\"oherer {O}rdnung.
   Abh. Math. Sem. Univ. Hamburg \textbf{36}, 140--149 (1971)
%
\bibitem{Frehse:1973:VI}
    Frehse, J.:
     On the regularity of the solution of the biharmonic
              variational inequality.
   Manuscripta Math. \textbf{9}, 91--103 (1973)
%
\bibitem{BGS:2018:POne}
  Brenner, S.C., Gedicke, J., Sung, L.-Y.:
  $P_1$ finite element methods for an elliptic optimal control problem with
                   pointwise state constraints.
  IMA J. Numer. Anal. (2018) doi:10.1093/imanum/dry071
%
\bibitem{Ladyzhenskaya:1958:Elliptic}
    Lady\v{z}enskaya, O.A.:
     On integral estimates, convergence, approximate methods, and
              solution in functionals for elliptic operators.
   Vestnik Leningrad. Univ.  \textbf{13}, 60--69 (1958)
%
\bibitem{BFS:1965:Element}
   Bogner, F.K.,  Fox, R.L.,  Schmit, L.A.:
  The generation of interelement compatible stiffness and
                   mass matrices by the use of interpolation formulas.
  In: Proceedings Conference on Matrix Methods in
                   Structural Mechanics, pp. 397--444.
  Wright Patterson A.F.B., Dayton, Ohio (1965)
%
\bibitem{AFS:1968:TUBA}
   Argyris, J.H., Fried, I., Scharpf, D.W.:
  The {TUBA} family of plate elements for the matrix
                  displacement method.
  Aero. J. Roy. Aero. Soc. \textbf{72}, 701--709 (1968)
%
\bibitem{CT:1965:Element}
  Clough, R.W., Tocher, J.L.:
  Finite element stiffbess matrices for analysis of plate bending.
    In: Proceedings Conference on Matrix Methods in
                   Structural Mechanics, pp. 515--545.
  Wright Patterson A.F.B., Dayton, Ohio (1965)
%
\bibitem{Ciarlet:1974:HCT}
  Ciarlet, P.G.:
  Sur l'\'el\'ement de Clough et Tocher.
  RAIRO Anal. Num\'er. \textbf{8}, 19--27 (1974)
%
\bibitem{DDPS:1979:Macro}
  Douglas J.Jr., Dupont, T., Percell, P., Scott, L.R.:
  A family of $C^1$ finite elements with optimal approximation
                   properties for various Galerkin methods for 2nd and 4th
                   order problems.
  R.A.I.R.O. Mod\'el. Math. Anal. Num\'er. \textbf{13}, 227--255 (1979)
%
@\bibitem{MB:1996:PTU}
    Melenk, J.M., Babu{\v{s}}ka, I.:
     The partition of unity finite element method: basic theory and
              applications
    Comput. Methods Appl. Mech. Engrg.
    \textbf{139}, 289--314 (1996)
%
\bibitem{BBO:2003:GFEM}
    Babu{\v{s}}ka, I. and Banerjee, U. and Osborn, J.E.:
    Survey of meshless and generalized finite element methods: a
              unified approach.
   Acta Numer. \textbf{12}, 1--125 (2003)
%
\bibitem{ODJ:2011:Plate}
  Oh,  H.S., Davis, C.B., Jeong, J.W.:
  Meshfree particle methods for thin plates.
  Comput. Methods Appl. Mech. Engrg. \textbf{209}, 156--171 (2012)
%
\bibitem{AC:1961:Plate}
   Adini, A., Clough, R.W.:
  Analysis of plate bending by the finite element method.
   NSF Report G. 7337 (1961)
%
\bibitem{BCIZ:1965:TEB}
   Bazeley, G.P.,  Cheung, Y.K., Irons, B.M., Zienkiewicz, O.C.:
   Triangular elements in bending
                    - conforming and nonconforming solutions.
                      In: Proceedings Conference on Matrix Methods in
                   Structural Mechanics, pp. 547--576.
  Wright Patterson A.F.B., Dayton, Ohio (1965)
%
\bibitem{Morley:1968:plate}
  Morley, L.S.D.:
  The triangular equilibrium problem in the solution of plate
                   bending problems.
  Aero. Quart. \textbf{19},  149--169 (1968)
%
\bibitem{FdV:1974:PatchTest}
  de Veubeke, B.F.:
  Variational principles and the patch test.
  Internat. J. Numer. Methods Engrg.
  \textbf{8}, 783--801 (1974)
%
\bibitem{Shi:1986:IBPlate}
   Shi, Z.-C.:
  On the convergence of the incomplete biquadratic
                    nonconforming plate element.
  Math. Numer. Sinica.  \textbf{8}, 53--62 (1986)
%
\bibitem{EGHLMT:2002:DG3D}
  Engel, G., Garikipati , K.,  Hughes, T.J.R.,
                   Larson, M.G.,  Mazzei, L. ,   Taylor, R.L.:
  Continuous/discontinuous finite element approximations
                   of fourth order elliptic problems in structural and
                   continuum mechanics with applications to thin beams
                   and plates, and strain gradient elasticity.
  Comput. Methods Appl. Mech. Engrg.
  \textbf{191}, 3669--3750 (2002)
%
\bibitem{BSung:2005:DG4}
   Brenner, S.C.,  Sung, L.-Y.:
   $C^0$ interior penalty methods for fourth order elliptic
                   boundary value problems on polygonal domains.
  J. Sci. Comput. \textbf{22/23}, 83--118 (2005)
%
\bibitem{Brenner:2012:C0IP}
  Brenner, S.C.:
  $C^0$ Interior Penalty Methods. In Blowey, J., Jensen, M. (eds.)
  Frontiers in Numerical Analysis-Durham 2010, pp. 79--147.
  Springer-Verlag, Berlin-Heidelberg (2012)
%
\bibitem{SM:2007:Biharmonic}
    S{\"u}li, E., Mozolevski, I.:
     $hp$-version interior penalty {DGFEM}s for the biharmonic
              equation.
   Comput. Methods Appl. Mech. Engrg. \textbf{196}, 1851--1863 (2007)
%
\bibitem{HHH:2010:C0DG}
    Huang, J., Huang, X., Han, W.:
    A new {$C^0$} discontinuous {G}alerkin method for
              {K}irchhoff plates.
   Comput. Methods Appl. Mech. Engrg.
   \textbf{199}, 1446--1454 (2010)
%
\bibitem{HH:2016:C0IP}
    Huang, X. and Huang, J.:
     A superconvergent {$C^0$} discontinuous {G}alerkin method for
              {K}irchhoff plates: error estimates, hybridization and
              postprocessing.
   J. Sci. Comput.  \textbf{69}, 1251--1278 (2016)
%
\bibitem{Meyer:2008:OptimalControl}
    Meyer, C.: Error estimates for the finite-element approximation of an
              elliptic control problem with pointwise state and control
              constraints. Control Cybernet. \textbf{37}, 51--83 (2008)
%
\bibitem{LGY:2009:Control}
    Liu, W., Gong, W., Yan, N.:
     A new finite element approximation of a state-constrained
              optimal control problem.
   J. Comput. Math.
    \textbf{27}, 97--114 (2009)
%
\bibitem{GY:2011:State}
  Gong, W., Yan, N.:
  A mixed finite element scheme for optimal control problems
                  with pointwise state constraints.
  J. Sci. Comput.
  \textbf{46}, 82--203 (2011)
 %
\bibitem{BSZ:2013:OptimalControl}
  Brenner, S.C., Sung, L.-Y. , Zhang, Y.:
  A quadratic {$C^0$} interior penalty method for an elliptic optimal control problem
                  with state constraints.
                  The IMA Volumes in Mathematics and its Applications. \textbf{157},
                  97--132 (2013)
 %
\bibitem{CMV:2014:State}
    Casas, E., Mateos, M., Vexler, B.:
     New regularity results and improved error estimates for
              optimal control problems with state constraints.
    ESAIM Control Optim. Calc. Var.
    \textbf{20}, 803--822 (2014)
%
\bibitem{BDS:2014:PUMOC}
  Brenner, S.C., Davis, C.B., Sung, L.-Y.:
  A partition of unity method for a class of fourth order elliptic variational
                       inequalities.
  Comp. Methods Appl. Mech. Engrg.
  \textbf{276}, 612--626 (2014)
%
\bibitem{BOPPSS:2016:OC3D}
  Brenner, S.C., Oh, M., Pollock, S., Porwal , K., Schedensack, M., Sharma, N.:
  A $C^0$ interior penalty method for elliptic distributed optimal control problems
                      in three dimensions with pointwise state constraints.
                      The IMA Volumes in Mathematics and its Applications.
                      \textbf{160}, 1--22 (2016)
%
\bibitem{Necas:2012:Direct}
 Ne{\v{c}}as, J.:Direct Methods in the Theory of Elliptic Equations,
 Springer, Heidelberg (2012)
%
\bibitem{GS:2002:Interpolation}
    Girault, V., Scott, L.R.:
     Hermite interpolation of nonsmooth functions preserving
              boundary conditions.
   Math. Comp. \textbf{71}, 1043--1074 (2002)
%
\bibitem{BSung:2017:State}
  Brenner, S.C., Sung, L.-Y.:
  A new convergence analysis of
                  finite element methods for elliptic distributed optimal
                   control problems with pointwise state constraints.
    SIAM J. Control Optim. \textbf{55}, 2289--2304 (2017)
%
\bibitem{Brenner:2015:CR}
 Brenner, S.C.:
  Forty years of the {C}rouzeix-{R}aviart element.
   Numer. Methods Partial Differential Equations.
   \textbf{31}, 367--396 (2015)
%
\bibitem{Wahlbin:2000:Handbook}
  Wahlbin, L.B.
  Local Behavior in Finite Element Methods.
   In: Ciarlet, P.G., Lions, J.L. (eds.)
  Handbook of Numerical Analysis, II, pp. 353--522
  North-Holland, Amsterdam (1991)
%
\bibitem{FGW:1973:SFM}
    Fix, G.J. , Gulati, S. , Wakoff, G.I.:
    On the use of singular functions with finite element
              approximations.
   J. Computational Phys. \textbf{13}, 209--228 (1973)
%
\bibitem{BKP:1979:DIE}
  Babu\v ska, I. , Kellogg, R.B., Pitk\"aranta, J.:
  Direct and inverse error estimates for finite
                   elements with mesh refinements.
   Numer. Math.  \textbf{33}, 447--471 (1979)
%
\bibitem{BSZ:2015:PP}
  Brenner, S.C.,  Sung, L.-Y. , Zhang, Y.:
  Post-processing procedures for a quadratic {${C^0}$} interior penalty
     method for elliptic distributed optimal control problems with pointwise state constraints.
  Appl. Numer. Math. \textbf{95}, 99--117  (2015)
%
\bibitem{BK:2002:PDAS}
    Bergounioux, M., Kunisch, K.:
    Primal-dual strategy for state-constrained optimal control
              problems.
   Comput. Optim. Appl.  \textbf{22}, 193--224 (2002)
%
\bibitem{HIK:2003:PDAS}
    Hinterm{\"u}ller, M., Ito, K., Kunisch, K.:
     The primal-dual active set strategy as a semismooth {N}ewton
              method.
   SIAM J. Optim.  \textbf{13}, 865--888 (2003)
%
\bibitem{IK:2008:Lagrange}
    Ito, K. and Kunisch, K.:
 Lagrange Multiplier Approach to Variational Problems and
              Applications.
  Society for Industrial and Applied Mathematics,
   Philadelphia (2008)
%
\bibitem{BM:2013:VEM}
  Brezzi, F., Marini, L.D.:
  Virtual element methods for plate bending problems.
 Comput. Methods Appl. Mech. Engrg. \textbf{253}, 455--462  (2013)
%
\bibitem{MWY:2014:WG4}
    Mu, L. and Wang, J. and Ye, X.:
     Weak {G}alerkin finite element methods for the biharmonic
              equation on polytopal meshes.
   Numer. Methods Partial Differential Equations.
    \textbf{30}, 1003--1029 (2014)
%
\bibitem{WW:2014:WG4}
    Wang, C. and Wang, J.:
     An efficient numerical scheme for the biharmonic equation by
              weak {G}alerkin finite element methods on polygonal or
              polyhedral meshes.
   Comput. Math. Appl. \textbf{68}, 2314--2330 (2014)
%
\bibitem{CM:2016:VEM}
    Chinosi, C., Marini, L.D.:
    Virtual element method for fourth order problems:
              {$L^2$}-estimates.
   Comput. Math. Appl. \textbf{72}, 1959--1967 (2016)
%
\bibitem{AMV:2018:VEM4}
    Antonietti, P.F. and Manzini, G. and Verani, M.:
    The fully nonconforming virtual element method for biharmonic
              problems.
   Math. Models Methods Appl. Sci. \textbf{28}, 387--407 (2018)
%
\bibitem{ZZCM:2018:MorleyVEM}
    Zhao, J. and Zhang, B. and Chen, S. and Mao, S.:
     The {M}orley-type virtual element for plate bending problems.
   J. Sci. Comput.  \textbf{76}, 610--629 (2018)
%
\bibitem{BDGK:2018:HHOPlate}
    Bonaldi, F., Di Pietro, D.A., Geymonat, G., Krasucki, F.:
    A hybrid high-order method for {K}irchhoff-{L}ove plate
              bending problems.
   ESAIM Math. Model. Numer. Anal.  \textbf{52}, 393--421 (2018)
%
\bibitem{BDR:2019:3DC1}
 Beir\~ao~da Veiga, L.,  Dassi, F., Russo,  A.:
 A {$C^1$} virtual element method on polyhedral meshes.
 arXiv:1808.01105v2 [math.NA] (2019)
%
\bibitem{BSZ:2019:Neumann}
    Brenner, S.C., Sung, L-Y., Zhang, Y.:
     $C^0$ interior penalty methods for an elliptic
              state-constrained optimal control problem with Neumann
              boundary condition.
   J. Comput. Appl. Math. \textbf{350}, 212--232 (2019)
%
\bibitem{BOS:2019:Neumann}
 Brenner, S.C., Oh, M., Sung, L.-Y.:
 $P_1$ finite element methods for an elliptic state-constrained distributed
 optimal control problem with Neumann boundary conditions.
 preprint (2019)
%
\bibitem{BGPS:2018:Morley}
   Brenner, S.C. , Gudi, T. and Porwal, K. and Sung, L.-Y. :
  A Morley finite element method for an elliptic distributed optimal control problem
                   with pointwise state and control constraints.
  ESAIM:COCV.  \textbf{24}, 1181--1206 (2018)
%
\bibitem{BGS:2018:Nonconvex}
     Brenner, S.C. , Gedicke, J., Sung, L.-Y.:
     $C^0$ interior penalty methods for an elliptic distributed
              optimal control problem on nonconvex polygonal domains with
              pointwise state constraints.
   SIAM J. Numer. Anal.  \textbf{56}, 1758--1785 (2018)
%
\bibitem{BWZ:2004:PF4}
  Brenner, S.C., Wang, K., Zhao, J.:
  Poincar\'e-Friedrichs inequalities for piecewise
                   $H^2$ functions.
  Numer. Funct. Anal. Optim. \textbf{25}, 463--478 (2004)
%
\bibitem{BNRS:2017:VonKarman}
  Brenner, S.C., Neilan, M., Reiser, A., Sung, L.-Y.:
  A $C^0$ interior penalty method for a von K\'arm\'an plate.
  Numer. Math. \textbf{135}, 803--832 (2017)
%
\bibitem{CC:2018:Obstacle}
    Christof, C. and Meyer, C.:
    A note on a priori {$L^p$}-error estimates for the obstacle
              problem.
   Numer. Math.  \textbf{139}, 27--45 (2018)
%
\bibitem{BGSZ:2017:AdaptiveVI4}
    Brenner, S.C., Gedicke, J., Sung, L.-Y.,  Zhang, Y.:
    An a posteriori analysis of $C^0$ interior penalty methods
              for the obstacle problem of clamped Kirchhoff plates.
   SIAM J. Numer. Anal.  \textbf{55}, 87--108 (2017)
%
%
\bibitem{BDS:2018:ASP}
    Brenner, S.C., Davis, C.B., Sung, L.-Y.:
  Additive {S}chwarz preconditioners for the obstacle problem of
              clamped {K}irchhoff plates.
   Electron. Trans. Numer. Anal. \text{49}, 274--290 (2018)
%
\bibitem{BDS:2018:DD25}
   Brenner, S.C., Davis, C.B., Sung, L.-Y.:
  Additive {Schwarz} preconditioners for a state
             constrained elliptic distributed optimal control
                    problem discretized by a partition of unity method.
  arXiv:1811.07809v1 [math.NA] (2018)
%
\end{thebibliography}


%
\end{document}